\documentclass[3p,preprint,review,12pt]{elsarticle}

\usepackage{amssymb}
\usepackage{mdwlist}
\usepackage[tbtags]{amsmath}
\usepackage{multirow}
\usepackage{graphics}
\usepackage{subfigure}
\usepackage{rotating}
\usepackage{enumitem}
\usepackage{epstopdf}
\usepackage{color}
\usepackage{epsf}

\newcommand{\Lagr}{\mathcal{L}}

\journal{Electric Power System Research}

\begin{document}

\begin{frontmatter}



\title{Impact of Equipment Failures and Wind Correlation on Generation Expansion Planning}

\author[dtu]{S.~Pineda\corref{cor1}} \ead{spinedamorente@gmail.com}
\author[dtu]{J.M.~Morales} \ead{juanmi82mg@gmail.com}
\author[dtu]{Y.~Ding} \ead{yding@elektro.dtu.dk}
\author[dtu]{J.~{\O}stergaard} \ead{joe@elektro.dtu.dk}
\cortext[cor1]{Corresponding author. Current address: University of Copenhagen, Universitetsparken 5, 2100 Copenhagen, Denmark. Phone: +45 45529800}
\address[dtu]{Technical University of Denmark, Elektrovej 325, 2800 Kgs. Lyngby, Denmark}


\begin{abstract}
Generation expansion planning has become a complex problem within a deregulated electricity market environment due to all the uncertainties affecting the profitability of a given investment. Current expansion models usually overlook some of these uncertainties in order to reduce the computational burden. In this paper, we raise a flag on the importance of both equipment failures (units and lines) and wind power correlation on generation expansion decisions. For this purpose, we use a bilevel stochastic optimization problem, which models the sequential and noncooperative game between the generating company (GENCO) and the system operator. The upper-level problem maximizes the GENCO's expected profit, while the lower-level problem simulates an hourly market-clearing procedure, through which LMPs are determined. The uncertainty pertaining to failures and wind power correlation are characterized by a scenario set, and their impact on generation expansion decisions are quantified and discussed for a 24-bus power system. 
\end{abstract}


\end{frontmatter}

%\linenumbers


\section{Notation}
\vspace{-2mm}
\subsection{Indexes and sets}
\begin{basedescript}{\desclabelstyle{\pushlabel}\desclabelwidth{2.5em}}
	\item[$b$] Index of load blocks.
	\item[$g/g'$] Index of existing/new conventional generating units.	
	\item[$n,m$] Index of buses.
	\item[$s$] Index of scenarios.
	\item[$w/w'$] Index of existing/new wind farms.	
	\item[$y$] Index of year of the planning horizon.
	\item[$\Gamma$] Set of existing wind farms owned by GENCO.		
  \item[$\Theta_{n}$] Set of wind farms connected to bus $n$.
	\item[$\Phi$] Set of existing conventional generating units owned by GENCO.	
	\item[$\Psi_{n}$] Set of conventional generating units connected to bus $n$.
	\item[$\Omega_{n}$] Set of buses connected to bus $n$.
\end{basedescript}	
\vspace{-3mm}
\subsection{Constants}
\begin{basedescript}{\desclabelstyle{\pushlabel}\desclabelwidth{3.4em}}
	\item[$B_{nm}$] Susceptance of line connecting buses $n$ and $m$ (p.u.).
	\item[$C^I_{g'/w'ny}$] Annualized investment cost of building unit $g'/w'$ at bus $n$ in year $y$ (\$).	
	\item[$C^P_{g/g'}$] Production cost of conventional generating unit $g/g'$ (\$/MWh).		
	\item[$k_{g/g's}$] Status of unit $g/g'$ in scenario $s$ (1 if available, 0 otherwise).
	\item[$k_{nms}$] Status of line $n-m$ in scenario $s$ (1 if available, 0 otherwise).
	\item[$L_b$] Load percentage (with respect to peak load) of block $b$ (p.u.).
	\item[$L_{nby}$] Load at bus $n$ corresponding to block $b$ of year $y$ (MW). 
	\item[$L^{peak}_{ny}$] Peak load at bus $n$ in year $y$ (MW). 	
	\item[$N^{T}_{w/w'}$] Number of wind turbines of wind farm $w/w'$.	
	\item[$T_b$] Duration of load block $b$ (h).
	\item[$\overline{P}^{C}_{g/g'}$] Capacity of conventional generating unit $g/g'$ (MW).
	\item[$\overline{P}^{F}_{nm}$] Capacity of line connecting buses $n$ and $m$ (MW).	
	\item[$\overline{P}^W_{ws}$] Power output of wind farm $w$ in scenario $s$ (MW).
	\item[$\overline{P}^{W,1}_{ws}$] Power output of a single turbine of wind farm $w$ in scenario $s$ (MW).
	\item[$\overline{P}^W_{w'ns}$] Power output of wind farm $w'$ in scenario $s$ if placed at bus $n$ (MW).
	\item[$\overline{P}^{W,1}_{w'ns}$] Power output of a single turbine within wind farm $w'$ in scenario $s$ if located at bus $n$ (MW).
	\item[$r$] Discount rate (\%).
	\item[$V^{L}$] Value of shed load (\$/MWh).
	\item[$\pi_{s}$] Probability of scenario $s$.	
\end{basedescript}	
\vspace{-3mm}
\subsection{Variables}
\begin{basedescript}{\desclabelstyle{\pushlabel}\desclabelwidth{3.0em}}
	\item[$L^{S}_{nsby}$] Amount of load shed from load block $b$ at bus $n$ in scenario $s$ in year $y$ (MW).	
	\item[$u_{g'/w'ny}$] Binary variable equal to 1 if unit $g'$/wind farm $w'$ is placed at bus $n$ in year $y$.		
	\item[$P^C_{g/g'sby}$] Dispatch of conventional unit $g/g'$ in scenario $s$ and load block $b$ in year $y$ (MW).		
	\item[$P^W_{w/w'sby}$] Dispatch of wind farm $w/w'$ in scenario $s$ and load block $b$ in year $y$ (MW).
	\item[$\delta_{nsby}$] Voltage angle at bus $n$, scenario $s$, and load block $b$ in year $y$ (rad).	
	\item[$\lambda_{nsby}$] Locational marginal price at bus $n$, scenario $s$, and load block $b$ in year $y$ (\$/MWh).	
	\item[$\Pi_{sby}$] GENCO's profit in scenario $s$ and load block $b$ in year $y$ (\$).	
\end{basedescript}

\section{Introduction}

The rapid growth of electricity demand in developed countries have turned the generation and transmission expansion planning into a key component in the long-term operation of power systems. In a \emph{vertically integrated} electricity supply industry, both generation and transmission expansion decisions are \emph{centrally undertaken to minimize the total cost}, including investment and operation cost, and/or maximize the reliability and security of the network. A review of the main methods for generation and transmission planning is presented in \cite{zhu1997review} and \cite{latorre2003classification}, respectively. Most of these models are formulated as one-level optimization problems involving a relatively low computational burden, even if uncertainties pertaining to demand level, wind power production, and equipment availability are accounted for \cite{kamalinia2011security,unsihuay2011multistage,Rocha201283, lopez2007generation,yu2009chance,Bakirtzis2012}.

On the other hand, since the \emph{liberalization} of the electricity sector in many countries around the world, new generation expansion models are required to determine GENCO's investment decisions according to \emph{profit maximization criteria} \cite{Baringo2013, botterud2007stochastic}. In this framework, GENCOs have to evaluate their expansion decisions according to the profits they make in the wholesale electricity market, which is cleared by an independent system operator (ISO) aiming, in turn, at minimizing the generation supply cost to satisfy electricity demand. Moreover, market outcomes are directly influenced by GENCOs' expansion decisions. It is easy to see that this problem falls within the framework of bilevel programming, in which two decision makers (GENCO and ISO), each with their individual objectives (maximize GENCO's profit and minimize system cost, respectively), act and react in a noncooperative sequential manner \cite{bard1998practical}. 

In this line, the bilevel model proposed in \cite{kazempourstrategic} calculates the optimal investment decisions for a producer to maximize its expected profit using a Benders decomposition approach. Similarly, reference \cite{wogrin2011generation} proposes a bilevel model that determines the expansion decisions of a producer via a conjectured-price response formulation. References \cite{careri2011generation} and \cite{baringo2012wind} analyze the generation capacity expansion including renewable energy sources (RES). In \cite{careri2011generation}, RES incentives such as feed-in tariff and quota obligation systems are considered. The stochastic bilevel model proposed in \cite{baringo2012wind} aims at identifying the new wind farms to be built to maximize the wind producer's profit, while risk considerations to determine expansion planning are discussed in \cite{Baringo2013}. Equilibrium models that analyze the generation expansion competitive behavior among GENCOs are discussed in \cite{murphy2005generation} and \cite{wang2009strategic}. Finally, reference \cite{roh2009market} proposes an iterative algorithm to analyze the coordination transmission and generation expansion made by the ISO and market agents.

Conversely to models that centrally determine the optimal expansion strategy, previously mentioned bilevel expansion models for GENCOs participating in electricity markets usually involve a significantly high computational complexity. In fact, although some of them include uncertainties corresponding to rival decisions \cite{wogrin2011generation,wang2009strategic}, demand variation \cite{kazempour2011strategic}, or wind power production \cite{baringo2011wind}, all of them disregard, to the best of our knowledge, the effect of equipment availability and wind correlation. 

In this paper, we propose a model to investigate how GENCO's optimal generation expansion decisions are influenced by the uncertainty corresponding to equipment failures and wind power correlation at different geographical locations. The proposed model consists of a bilevel formulation in which the upper-level problem maximizes the GENCO's expected profit including both pool revenue and investment cost. Pool revenue is computed according to the locational marginal prices (LMP) and dispatched quantities resulting from the lower-level market-clearing problem which, in turn, aims at minimizing expected operational cost while ensuring the continuous balance between production and consumption and the fulfillment of the network constraints. Using the primal-dual theory, a single-level mixed-integer equivalent formulation that can be readily solved by off-the-shelf optimization software is obtained.

The main contributions of this paper are thus twofold: 

\begin{enumerate}[leftmargin=0.4cm,itemindent=0cm]
\item To investigate the impact of unit and line failures on optimal generation expansion decisions of conventional generating units.
\item To investigate the impact of wind power correlation on optimal generation expansion decisions of wind farms.
\end{enumerate}

This paper is organized as follows. Section \ref{SectionModel} elaborates on the model assumptions. The bilevel optimization problem and its equivalent MILP formulation are described in Section \ref{SectionFormulation}. In Section \ref{SectionCaseStudy}, the results of a 24-bus case study are provided. The main conclusions of this work are discussed in Section \ref{SectionConclusions}. Finally, \ref{Appendix} contains linearization technicalities.

\section{Model assumptions}\label{SectionModel}

Transmission and generation expansion planning can be classified as dynamic or static according to the treatment of the planning horizon \cite{akbari2012security}. Dynamic expansion planning includes several years and investment decisions are sequentially made along the decision horizon. Parameter changes throughout the planning period as well as the use of discount rate to compute the net present value of future costs can be readily incorporated in a dynamic approach \cite{kamalinia2011security,unsihuay2011multistage, Bakirtzis2012, Baringo2013, botterud2007stochastic,wogrin2011generation, careri2011generation,roh2009market,akbari2012security}. Dynamic expansion models demand, though, a high computational effort, thus requiring some simplifications to keep them computationally tractable \cite{latorre2003classification} or the use of decomposition techniques \cite{Baringo2013}. Alternatively to dynamic expansion planning, static models are computationally more manageable because expansion decisions are determined for a single future target year. One-year static expansion models have been used for generation \cite{kazempourstrategic, baringo2012wind,wang2009strategic,chung2004optimal}, transmission \cite{sanchez2005probabilistic,garces2009bilevel}, as well as combined generation and transmission planning \cite{lopez2007generation}. 

In this paper, a multi-year formulation to determine the optimal expansion decisions of a power producer taking into account the impact of equipment failures and wind correlation is proposed. In some analysis though, expansion decisions are investigated using a single target year in order to maintain the computational burden of the model within reasonable limits.
%However, since this paper primary focuses on quantifying the effect of equipment failures and wind speed correlation on GENCO's expansion planning, a single future target year is considered in the case study in order to keep the computational time within reasonable levels. 

In order to take advantage of economies of scale, typical constructions of new units are given in discrete and relatively large amounts. We model this fact by just considering a finite set of possible new investments and modeling expansion decisions through the binary variables $u_{g'ny}$ and $u_{w'ny}$, which are equal to 1 if unit $g'$ and wind farm $w'$, respectively, are built at bus $n$ in year $y$, and 0 otherwise. The binary variables $\hat{u}_{g'ny}$ and $\hat{u}_{w'ny}$ are equal to 1 if unit $g'$ or wind farm $w'$ are built at or before year $y$, and 0 otherwise. $\overline{P}^{C}_{g'}$ and $C^P_{g'}$ represent the capacity and production cost of new conventional units. Likewise, new wind farms are only characterized by their capacity $\overline{P}^{W}_{w'}$, being their marginal cost equal to 0. The investment cost of new units and wind farms for each year $y$ are denoted by $C^I_{g'ny}$ and $C^I_{w'ny}$, in that order. 

Assuming a perfectly competitive market, selling offers submitted by producers exactly represent their corresponding marginal costs \cite{akbari2012security}. Moreover, for the sake of simplicity, we suppose that all the revenues produced by new and old units come solely from the sale of electrical energy in the wholesale market, then disregarding additional revenues associated with ancillary services or capacity payments. 

Since uncertainty pertaining to demand or bidding strategies are  modelled in detail in the technical literature, and in order to render the subsequent analysis and discussion more intuitive, the modeling of uncertainties different from equipment availability and wind speed have not been considered here. Note, however, that such uncertainties can be easily included in the proposed model by just increasing the number of scenarios.

Wind speed stochastic behavior is modeled using historical data at different geographical locations. Once the wind speed scenarios are selected, the power production of a single wind turbine ($\overline{P}^{W,1}_{ws}$ and $\overline{P}^{W,1}_{w'ns}$) is computed according to its power curve. Thus, the power production of existing and new wind farms is determined as 
\begin{equation}
 \overline{P}^W_{ws} = N^{T}_{w}\overline{P}^{W,1}_{ws}; \quad\quad\quad   \overline{P}^W_{w'ns} = N^{T}_{w'}\overline{P}^{W,1}_{w'ns}. 	
\end{equation}

The availability status of existing/new units and transmission lines are modeled through parameters $k_{gs}$, $k_{g's}$, and $k_{nms}$, respectively, being equal to 1 if the corresponding device is available in scenario $s$, and 0 otherwise. In order to keep the model computationally tractable, only single-element contingencies are considered \cite{akbari2012security,billinton1996reliability}. As an example, the probability that unit $g_1$ fails is determined as
\begin{equation}
\hspace{0mm}\textsc{FOR}_{g_1}\cdot\hspace{-1mm}\prod_{g \neq g_1}(1-\textsc{FOR}_g)\cdot\prod_{g'}(1-\textsc{FOR}_{g'})\cdot\hspace{-4mm} \prod_{nm:m\in \Omega_{n}}\hspace{-4mm}(1-\textsc{FOR}_{nm})	\label{ProbConting}
\end{equation}

\noindent where $\textsc{FOR}_{g/g'/nm}$ represents the probability of having an individual unexpected failure of device $g/g'/nm$. Since not all the possible contingencies are considered, the probability calculated in \eqref{ProbConting} has to be normalized. Note, however, that including multi-device failure scenarios in the analysis is straightforward. For simplicity, unexpected failures of wind farms are not considered in this analysis. 

\begin{figure}	\centering	
		\includegraphics[scale=0.5]{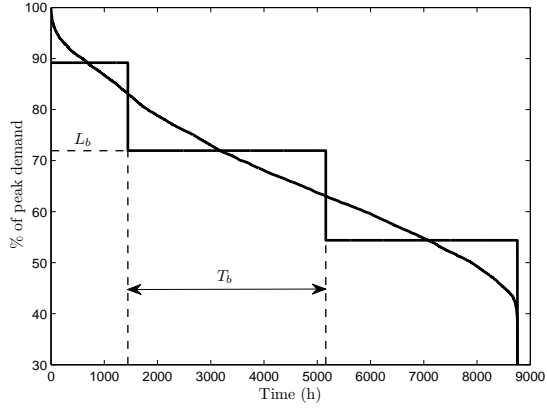}\vspace{-3mm} 		
	\caption{Yearly load duration curve approximation}\vspace{-3mm}	
	\label{fig:LoadDemandCurve}	
\end{figure}

The lower-level optimization problem consists of an hourly market-clearing algorithm that minimizes the operational cost including DC power flow equations as well as maximum capacity constraints for both generating units and transmission lines. The system reliability cost is modeled through the value of shed load $V^{L}$. Electricity consumption is considered to be known, inelastic, and uncorrelated with the wind power production. The hourly load duration curve is approximated by $N_B$ blocks (as depicted in Fig.\,\ref{fig:LoadDemandCurve}) in order to account for the load variability throughout each year of the planning horizon \cite{Bakirtzis2012}. Assuming that the load is proportionally distributed among all buses, the load at each bus $n$, block $b$ and year $y$ is computed as
\begin{equation}
	L_{nby} = L_b \cdot L^{peak}_{ny}, \quad \forall n, \forall b, \forall y,
\end{equation}

\noindent where $L^{peak}_{ny}$ is the peak load at each bus of the network in year $y$.

\section{Model Formulation}\label{SectionFormulation}

The multi-year bilevel stochastic optimization problem used to determine the impact of equipment failures and wind correlation on the GENCO's expansion planning decisions is presented below. 
\begin{subequations} \label{CapExp1}
\begin{gather}
\textrm{Maximize}_{u_{g'ny},u_{w'ny}} \nonumber\\
\sum_y \frac{1}{(1+r)^y} \Bigg( \sum_{sb}\pi_{s}T_b\Pi_{sby} - \sum_{n}\Big( \sum_{g'}\hat{u}_{g'ny}C^I_{g'ny} + \sum_{w'}\hat{u}_{w'ny}C^I_{w'ny} \Big) \Bigg)  \label{CapExp1_OF}\displaybreak[1]
\end{gather}
subject to
\begin{align}
& \sum_{ny} u_{g'ny} \leq 1, \quad  \forall  g' \label{CapExp1_OneUnit}\displaybreak[1]\\
& \sum_{ny} u_{w'ny} \leq 1, \quad  \forall  w' \label{CapExp1_OneWindFarm}\displaybreak[1]\\
& \hat{u}_{g'ny} = \sum_{y' \leq y}u_{g'ny'}, \quad \forall g', \forall n, \forall y \label{CapExp1_CumCapUnit}\displaybreak[1]\\
& \hat{u}_{w'ny} = \sum_{y' \leq y}u_{w'ny'}, \quad \forall w', \forall n, \forall y \label{CapExp1_CumCapWind}\displaybreak[1]\\
& \Pi_{sby} = \sum_{g \in \Phi}(\lambda_{nsby:g\in\Psi_{n}}-C^P_{g})P^C_{gsby}  + \sum_{w \in \Gamma} \lambda_{nsby:w\in\Theta_{n}}P^W_{wsby}+ \nonumber\\
& + \sum_{g'n} (\hat{u}_{g'ny}\lambda_{nsby}-C^P_{g'})P^C_{g'sby} + \sum_{w'n} \hat{u}_{w'ny}\lambda_{nsby}P^W_{w'sby}, \quad \forall  s, \forall b, \forall y \label{CapExp1_ProfitGenCo}\displaybreak[1]\\ 
& (P^C_{g/g'sby},P^W_{w/w'sby},\lambda_{nsby}) \in  \textrm{arg} \Bigg\{ \textrm{Min}  \sum_{g}C^P_{g}P^C_{gsby} +\sum_{g'}C^P_{g'}P^C_{g'sby} + \sum_{n}V^{L} L^{S}_{nsby} \label{CapExp1_SocialWelfare}\displaybreak[1]\\
& \textrm{subject to} \nonumber\\
&  0 \leq P^C_{gsby} \leq k_{gs}\overline{P}^C_{g}:\phi^{min}_{gsby},\phi^{max}_{gsby}, \quad \forall g \label{CapExp1_MaxPG}\displaybreak[1]\\
&  0 \leq P^C_{g'sby} \leq k_{g's}\sum_n\hat{u}_{g'ny}\overline{P}^{C}_{g'}:\phi^{min}_{g'sby},\phi^{max}_{g'sby}, \quad \forall g' \label{CapExp1_MaxPg'}\displaybreak[1]\\
& 0 \leq L^{S}_{nsby} \leq L_{nby} : \beta^{min}_{nsby},\beta^{max}_{nsby}, \quad \forall n \label{CapExp1_MaxShed}\displaybreak[1]\\
& 0 \leq P^W_{wsby} \leq  \overline{P}^W_{ws}: \gamma^{min}_{wsby},\gamma^{max}_{wsby}, \quad \forall w \label{CapExp1_MaxSpill}\displaybreak[1]\\
& 0 \leq P^W_{w'sby} \leq  \sum_{n}\hat{u}_{w'ny}\overline{P}^W_{w'ns}  : \gamma^{min}_{w'sby},\gamma^{max}_{w'sby}, \quad \forall w'  \label{CapExp1_MaxSpill2}\displaybreak[1]\\
& \sum_{g\in\Psi_{n}}P^C_{gsby}  + \sum_{g'}\hat{u}_{g'ny}P^C_{g'sby}  +  \sum_{w \in \Theta_{n}}P^{W}_{wsby}  + \sum_{w'}\hat{u}_{w'ny}P^{W}_{w'sby}  = \nonumber\\
& \qquad = L_{nby} - L^{S}_{nsby} +  \sum_{m\in\Omega_{n}}B_{nm}k_{nms}(\delta_{nsby}-\delta_{msby}):\lambda_{nsby},  \quad \forall  n \label{CapExp1_Balance}\displaybreak[1]\\
& B_{nm}k_{nms}(\delta_{nsby}-\delta_{msby}) \leq k_{nms}\overline{P}^{F}_{nm}:\theta^{max}_{nmsby},  \forall  n,m   \label{CapExp1_MaxFlow}\displaybreak[1] \\
& \delta_{n_1sby}=0: \xi_{n_1sby} \label{CapExp1_SlackBus} \displaybreak[1] \Bigg\}, \quad \forall s, \forall b, \forall y.
\end{align}
\end{subequations}\vspace{-0mm}

Model \eqref{CapExp1_OF}--\eqref{CapExp1_SlackBus} is a bilevel optimization problem. The upper-level objective function \eqref{CapExp1_OF} aims at maximizing the discounted GENCO's expected profit throughout a multi-year planning horizon subject to a set of lower-level optimization problems, \eqref{CapExp1_SocialWelfare}--\eqref{CapExp1_SlackBus}, representing an hourly pool-based market clearing, one for each scenario $s$, load block $b$ and year $y$. Objective function \eqref{CapExp1_OF} includes the sum of the discounted profit from selling electricity in the pool (first term) and the investment cost of new units and wind farms (second term). Equations \eqref{CapExp1_OneUnit} and \eqref{CapExp1_OneWindFarm} impose that each unit and wind farm can be either placed at one single bus in a specific year or not built at all. Equations \eqref{CapExp1_CumCapUnit} and \eqref{CapExp1_CumCapWind} ensure that $\hat{u}_{g'ny}$ and $\hat{u}_{w'ny}$ are equal to 1 if a new unit or wind farm is built at of before year $y$, respectively. Equation \eqref{CapExp1_ProfitGenCo} computes the pool profit for each scenario $s$ load block $b$ and year $y$ as the sum of the profit corresponding to existing units (first term), existing wind farms (second term), new units (third term), and new wind farms (fourth term). 

The lower-level objective function \eqref{CapExp1_SocialWelfare} minimizes the system cost for each scenario $s$, load block $b$ and year $y$. Constraints \eqref{CapExp1_MaxPG} and \eqref{CapExp1_MaxPg'} limit the output of each unit. Constraints  \eqref{CapExp1_MaxShed}, \eqref{CapExp1_MaxSpill} and \eqref{CapExp1_MaxSpill2} limit, respectively, the load shedding at each bus and the wind dispatch of existing and new wind farms. Constraint \eqref{CapExp1_Balance} is the power balance equation at each bus. Equation \eqref{CapExp1_MaxFlow} imposes the maximum power flow through the  lines. Equation \eqref{CapExp1_SlackBus} arbitrarily sets the value of angle $\delta_{n_1}$ to 0. Note that dual variables corresponding to the constraints of the lower-level problems are included after a colon.

Note that although not considered here for the sake of simplicity, the multi-year formulation \eqref{CapExp1_OF}--\eqref{CapExp1_SlackBus} can also be modified to incorporate additional features such as constructions times of the generating units or risk measures.

To solve the bilevel optimization problem \eqref{CapExp1_OF}--\eqref{CapExp1_SlackBus} we replace each lower-level problem by its corresponding primal and dual constraints plus the strong duality theorem \cite{garces2009bilevel,motto2005mixed}.
\begin{subequations} \label{CapExp2}
\begin{gather}
\textrm{Maximize}_{u_{g'ny},u_{w'ny}} \nonumber\\
\sum_y \frac{1}{(1+r)^y} \Bigg( \sum_{sb}\pi_{s}T_b\Pi_{sby} - \sum_{n}\Big( \sum_{g'}\hat{u}_{g'ny}C^I_{g'ny} + \sum_{w'}\hat{u}_{w'ny}C^I_{w'ny} \Big) \Bigg)  \label{CapExp2_OF}\displaybreak[1]
\end{gather}
subject to
\begin{align}
& \eqref{CapExp1_OneUnit}-\eqref{CapExp1_ProfitGenCo} \nonumber\displaybreak[1]\\
& \eqref{CapExp1_MaxPG}-\eqref{CapExp1_SlackBus} \nonumber\displaybreak[1]\\
& \lambda_{nsby:g \in\Psi_{n}} + \phi^{max}_{gsby} + \phi^{min}_{gsby} = C^P_{g}, \quad \forall g, \forall s, \forall b, \forall y \label{CapExp2_Dual1}\displaybreak[1]\\
& \sum_{n}\hat{u}_{g'ny}\lambda_{nsby} + \phi^{max}_{g'sby} + \phi^{min}_{g'sby} = C^P_{g'}, \quad \forall g', \forall s, \forall b, \forall y \label{CapExp2_Dual2}\displaybreak[1]\\
& \lambda_{nsby} + \beta^{max}_{nsby} + \beta^{min}_{nsby} = V^{L}, \quad  \forall  n, \forall s, \forall b, \forall y \label{CapExp2_Dual3}\displaybreak[1]\\ 
& -\lambda_{nsby:w \in \Theta_{n}} + \gamma^{max}_{wsby} + \gamma^{min}_{wsby} = 0, \quad  \forall  w, \forall s, \forall b, \forall y \label{CapExp2_Dual4}\displaybreak[1]\\
& -\sum_{n}\hat{u}_{w'ny}\lambda_{nsby} + \gamma^{max}_{w'sby} + \gamma^{min}_{w'sby} = 0, \quad  \forall  w', \forall s, \forall b, \forall y \label{CapExp2_Dual4_2}\displaybreak[1]\\
& \sum_{m \in \Omega_{n}}  B_{nm} (\lambda_{msby}-\lambda_{nsby}+\theta^{max}_{nmsby} - \theta^{max}_{mnsby})+\xi_{n_1sby} =0, \quad \forall n, \forall s, \forall b, \forall y \label{CapExp2_Dual5}\displaybreak[1]\\
& \phi^{max}_{gsby}, \phi^{max}_{g'sby}, \beta^{max}_{nsby}, \theta^{max}_{mnsby}, \gamma^{max}_{wsby}, \gamma^{max}_{w'sby} \leq 0 \label{CapExp2_Dual6}\displaybreak[1]\\ 
& \phi^{min}_{gsby}, \phi^{min}_{g'sby}, \beta^{min}_{nsby}, \gamma^{min}_{wsby}, \gamma^{min}_{w'sby} \geq 0 \label{CapExp2_Dual7}\displaybreak[1]\\
& \sum_{g}C^P_{g}P^C_{gsby} + \sum_{g'}C^P_{g'}P^C_{g'sby} +  \sum_{n}V^{L} L^{S}_{nsby} =  \sum_{g}\phi^{max}_{gsby} k_{gs}\overline{P}^{C}_{g} + \nonumber\\ 
& +\sum_{g'n}\phi^{max}_{g'sby} \hat{u}_{g'ny}k_{g's}\overline{P}^{C}_{g'}+ \sum_{n}\beta^{max}_{nsby}L_{nby} + \sum_{w} \gamma^{max}_{wsby}\overline{P}^W_{ws}  + \nonumber\\
& +  \sum_{w'n}\gamma^{max}_{w'sby}\hat{u}_{w'ny}\overline{P}^W_{w'ns} + \sum_{n,m\in \Omega_{n}}k_{nms}\overline{P}^F_{nm}\theta^{max}_{nmsby}+  \sum_{n}\lambda_{nsby}L_{nby}, \forall s, \forall b, \forall y. \label{CapExp2_StrongDuality}\displaybreak[1]
\end{align}
\end{subequations}

Equations \eqref{CapExp1_OneUnit}--\eqref{CapExp1_ProfitGenCo} are the constraints corresponding to the upper-level problem. Equations \eqref{CapExp1_MaxPG}--\eqref{CapExp1_SlackBus} represent the primal constraints of the lower-level problem. Likewise, equations \eqref{CapExp2_Dual1}--\eqref{CapExp2_Dual7} correspond to the constraints of the dual formulation of the lower-level problem. Finally, constraints \eqref{CapExp2_StrongDuality} ensure that the primal and dual formulation of the lower-level problems reach the same objective function at the optimal solution. The formulation above contains several non-linear terms:
\begin{enumerate}[leftmargin=0.4cm,itemindent=0cm]
\item Product of binary variables and continuous variables: $\hat{u}_{g'ny}\lambda_{nsby}$ in equation \eqref{CapExp2_Dual2}, $\hat{u}_{g'ny}P^C_{g'sby}$ and $\hat{u}_{w'ny}P^W_{w'sby}$ in equation \eqref{CapExp1_Balance}, $\hat{u}_{w'ny}\lambda_{nsby}$ in equation \eqref{CapExp2_Dual4_2}, and $\hat{u}_{g'ny}\phi^{max}_{g'sby}$ and $\hat{u}_{w'ny}\gamma^{max}_{nsby}$ in equation \eqref{CapExp2_StrongDuality}. Note that these expressions can be linearized as explained in \ref{Appendix} \cite{floudas1995nonlinear}.
\item Product of two continuous variables: $\lambda_{nsby:g\in\Psi_{n}}P^C_{gsby}$, $\lambda_{nsby:w\in\Theta_{n}}P^W_{wsby}$, $\hat{u}_{g'ny}\lambda_{nsby}$ $P^C_{g'sby}$ and $\hat{u}_{w'ny}\lambda_{nsby}P^W_{w'sby}$, all of them in equation \eqref{CapExp1_ProfitGenCo}. The procedure to linearize these terms using KKT conditions is explained in \ref{Appendix}.
\end{enumerate}

In doing so, optimization model \eqref{CapExp2_OF}--\eqref{CapExp2_StrongDuality} can be equivalently formulated as a mixed-integer linear optimization problem that can be solved using commercial software. 

\section{Case study}\label{SectionCaseStudy}

In this section, we use formulation \eqref{CapExp2_OF}--\eqref{CapExp2_StrongDuality} to determine the optimal generation expansion decisions of a power producer considering both a single target year and a planning horizon of three years.

\subsection{Data}

The IEEE RTS 24-bus system is analyzed in this section \cite{grigg1999ieee}. The data corresponding to already existing units and transmission lines are provided in Table \ref{table:unitdata} and Table \ref{table:linedata}, respectively. Table \ref{table:peakdata} provides the peak load at each bus. The total system load throughout one year is approximated by the 20 blocks provided in Table \ref{table:loadblocks}. It is worth mentioning that under any single-device failure scenario, no load shedding occurs in the system, i.e., the power system satisfies the N-1 reliability criterion. 

{\setlength{\tabcolsep}{1.2mm}\begin{table}[htb]\begin{center}\renewcommand{\arraystretch}{1}
\caption{Generating unit data}\label{table:unitdata}
\begin{tabular}{c c c c c || c c c c c }
  \hline
  $n$    & $\# units$ & $\overline{P}^C_{g}$ & $C^P_{g}$ & $\textrm{FOR}_g$ & $n$    & $\# units$ & $\overline{P}^C_{g}$ & $C^P_{g}$ & $\textrm{FOR}_g$\\
  \hline
  $n_1$    & 2					& 20						& 43.5		& 10 & $n_{15}$ & 1					& 155						& 11.5		& 4\\
  $n_1$    & 2					& 76						& 14.4		& 2  & $n_{16}$ & 1					& 155						& 11.5		& 4\\
  $n_2$    & 2					& 20						& 43.5		& 10 & $n_{18}$ & 1					& 400						& 6.0 		& 12 \\
  $n_2$    & 2					& 76						& 14.4		& 2  & $n_{21}$ & 1					& 400						& 6.0 		& 12\\
  $n_7$    & 3					& 100						& 23.0		& 4  & $n_{22}$ & 6					& 50						& 0.0 		& 1 \\
  $n_{13}$ & 3					& 197						& 22.1		& 5  & $n_{23}$ & 2					& 155						& 11.5		& 4\\
  $n_{15}$ & 5					& 12						& 27.6		& 2  & $n_{23}$ & 1					& 350					  & 11.4		& 8\\  
  \hline
\end{tabular}\vspace{-3mm}\end{center}
\end{table}
}
{\setlength{\tabcolsep}{0.7mm}\begin{table}[htb]\begin{center}\renewcommand{\arraystretch}{1}
\caption{Network data}\label{table:linedata}
\begin{tabular}{c c c || c c c || c c c || c c c}
\hline
$nm$ 	& $B_{nm}$  & $\overline{P}^F_{nm}$  & $nm$ 	& $B_{nm}$  & $\overline{P}^F_{nm}$  & $nm$   & $B_{nm}$  & $\overline{P}^F_{nm}$ & $nm$   & $B_{nm}$  & $\overline{P}^F_{nm}$ \\
  \hline  
  $n_{1}n_{2}$	  & 68.5 & 50  & 	$n_{6}n_{10}$	  & 15.6  & 150  & $n_{11}n_{14}$  & 23.5	& 200 & $n_{16}n_{17}$ & 38.0	& 200\\
  $n_{1}n_{3}$	  & 4.4  & 75  & 	$n_{7}n_{8}$	  & 15.3  & 200  & $n_{12}n_{13}$  & 20.5	& 225 & $n_{16}n_{19}$ & 42.7	& 300\\
  $n_{1}n_{5}$	  & 11.0 & 100 & 	$n_{8}n_{9}$	  & 5.7   & 100  & $n_{13}n_{23}$  & 10.2	& 200 & $n_{17}n_{18}$ & 69.9	& 125\\
  $n_{2}n_{4}$	  & 7.4  & 75  & 	$n_{8}n_{10}$   & 5.7	  & 100  & $n_{13}n_{23}$  & 11.3	& 125 & $n_{17}n_{22}$ & 9.4	  & 75\\
  $n_{2}n_{6}$	  & 4.9  & 150 & 	$n_{9}n_{11}$   & 11.9	& 150  & $n_{14}n_{16}$  & 16.8	& 250 & $n_{18}n_{21}$ & 75.8	& 275\\
  $n_{3}n_{9}$	  & 7.9  & 125 & 	$n_{9}n_{12}$   & 11.9	& 175  & $n_{15}n_{16}$  & 58.1	& 200 & $n_{19}n_{20}$ & 49.3	& 275\\
  $n_{3}n_{24}$	  & 11.9 & 175 & 	$n_{10}n_{11}$  & 11.9	& 200  & $n_{15}n_{21}$  & 40.2	& 275 & $n_{20}n_{23}$ & 89.3	& 400\\
  $n_{4}n_{9}$	  & 9.0  & 75  & 	$n_{10}n_{12}$  & 11.9	& 225  & $n_{15}n_{24}$  & 18.9	& 175 & $n_{21}n_{22}$ & 14.5	& 75\\
  $n_{5}n_{10}$	  & 10.6 & 75  & 	$n_{11}n_{13}$  & 20.5	& 325  &  &  & 	& 	& &  \\ 
  \hline
\end{tabular}\vspace{-3mm}\end{center}
\end{table}}

{\setlength{\tabcolsep}{1mm}\begin{table}[htb]\begin{center}\renewcommand{\arraystretch}{1}
\caption{Peak load data}\label{table:peakdata}
\begin{tabular}{c c c c c c c c c c c c c c c c c c}
  \hline
  $n$ & $n_{1}$ & $n_{2}$& $n_{3}$& $n_{4}$& $n_{5}$& $n_{6}$& $n_{7}$ & $n_{8}$ & $n_{9}$ & $n_{10}$ & $n_{13}$& $n_{14}$& $n_{15}$& $n_{16}$& $n_{18}$& $n_{19}$& $n_{20}$ \\
  $L^{peak}_{n}$ & 108 & 97 & 180 & 74 & 71 & 137 & 125 & 171 & 174 & 194 & 265 & 194 & 316 & 100 & 334 & 182 & 128 \\  
  \hline
\end{tabular}\vspace{-3mm}\end{center}
\end{table}
}

{\setlength{\tabcolsep}{2.5mm}\begin{table}[htb]\begin{center}\renewcommand{\arraystretch}{1}
\caption{Load blocks}\label{table:loadblocks}
\begin{tabular}{c c c c c c c c c c c}
  \hline
  $b$ & $b_{1}$ & $b_{2}$& $b_{3}$& $b_{4}$& $b_{5}$& $b_{6}$& $b_{7}$ & $b_{8}$ & $b_{9}$ & $b_{10}$ \\
  $L_{b}$ & 0.47 & 0.49 & 0.52 & 0.55 & 0.58 & 0.61 & 0.63 & 0.66 & 0.69 & 0.72 \\ 
  $T_{b}$ & 52 & 291 & 287 & 570 & 647 & 511 & 379 & 530 & 708 & 543 \\ 
  \hline
  $b$ & $b_{11}$& $b_{12}$& $b_{13}$& $b_{14}$& $b_{15}$& $b_{16}$& $b_{17}$ & $b_{18}$ & $b_{19}$ & $b_{20}$\\
  $L_{b}$ & 0.75 & 0.77 & 0.80 & 0.83 & 0.86 & 0.89 & 0.92 & 0.94 & 0.97 & 1.00 \\  
  $T_{b}$ & 196 & 147 & 183 & 258 & 485 & 751 & 760 & 823 & 456 & 183 \\
  \hline
\end{tabular}\vspace{-3mm}\end{center}
\end{table}
}

\subsection{Impact of unit and line failures on capacity expansion}\label{sec:ImpFailures}

The impact of unit and line failures on the expansion of conventional units is analyzed and discussed next. For simplicity, no wind farms are included in the network. Unavailability rates of generating units are provided in Table \ref{table:unitdata}, while the unavailability rate of all transmission lines is considered equal to 2\%. Assuming that new units can be built at any bus of the transmission grid, the effect of failures on GENCO's profit is computed as follows:
\begin{enumerate}
	\item Optimization problem \eqref{CapExp2_OF}--\eqref{CapExp2_StrongDuality} is solved considering that all units and lines are available, denoting the locations of the new units as $B^{\text{\tiny NF}}$.
	\item GENCO's expected profit if expansion decisions are made without considering failures, $\overline{\Pi}^{\text{\tiny NF}}$, is determined by solving \eqref{CapExp2_OF}--\eqref{CapExp2_StrongDuality} with the complete availability scenario set and the new unit locations fixed to $B^{\text{\tiny NF}}$.
	\item Optimization problem \eqref{CapExp2_OF}--\eqref{CapExp2_StrongDuality} is solved again including all the availability scenarios. The optimal location and expected profit are denoted by $B^{\text{\tiny F}}$ and $\overline{\Pi}^{\text{\tiny F}}$, respectively. 
  \item The impact of failures on expansion decisions is measured through the difference between $\overline{\Pi}^{\text{\tiny NF}}$ and $\overline{\Pi}^{\text{\tiny F}}$, which is denoted as $\Delta\overline{\Pi}$.
\end{enumerate}

The above procedure is carried out to determine the impact of equipment failures on the expansion decisions of a GENCO that can decide on the installation of four generating units of 50 MW each and the same marginal cost $C^P_{g'}$. Table \ref{table:impactfailure} presents the results for different values of the marginal cost $C^P_{g'}$ considering a single target year. For simplicity, the expected profit (expressed in \$million) is computed assuming that the GENCO only owns the new unit and the annualized investment cost is proportional to its capacity (\$400/kW to be paid in 40 years) and independent of its location. New generating units can be installed in buses $n_1$, $n_2$ and $n_7$.

{\setlength{\tabcolsep}{3mm}\begin{table}[htb]\begin{center}\renewcommand{\arraystretch}{1}
\caption{Impact of failures on expansion planning (case study)}\label{table:impactfailure}
\begin{tabular}{c c c c c c c c}
  \hline 
   $C^P_{g'}$ & $B^{\text{\tiny NF}}$ & $\overline{\Pi}^{\text{\tiny NF}}$ & $B^{\text{\tiny F}}$ & $\overline{\Pi}^{\text{\tiny F}}$ & $\Delta\overline{\Pi}$ \\
   \hline
15  & $n_2,n_2,n_7,n_7$ & 3.93 & $n_1,n_7,n_7,n_7$ 	& 4.42 & 12.3\%     \\
16  & $n_1,n_2,n_7$  	& 2.61 & $n_1,n_7,n_7,n_7$ 	& 3.04 & 16.5\%     \\
17  & $n_2,n_7$  		& 1.72 & $n_2,n_7$   		& 1.72 & 0\%      \\
18  & $n_7$  			& 0.63 & $n_1,n_7$			& 0.97 & 53.1\%     \\
19  & $n_7$  			& 0.31 & $n_1,n_7$			& 0.43 & 38.2\%     \\
20  & $n_7$  			& 0.10 & $n_2$ 				& 0.27 & 159.6\%     \\
21  & $-$  				& $-$  & $n_2$ 				& 0.19 & $-$     \\
22  & $-$  				& $-$  & $n_2$ 				& 0.11 & $-$     \\
23  & $-$  				& $-$  & $n_1$ 				& 0.04 & $-$     \\
24  & $-$  				& $-$  & $n_2$ 				& 0.01 & $-$     \\
  \hline
\end{tabular}\vspace{-3mm}\end{center}
\end{table}
}

Note that in spite of the low outage rates of generating units and transmission lines, the expected profit of a GENCO deciding the optimal location of a generating units of 50MW with a marginal cost of 18\$/MWh increases 53\% if unexpected failures are considered in the generation expansion problem. Furthermore, observe that the modeling of equipment failures makes profitable the installation of generating units with a marginal cost between 21\$/MWh and 24\$/MWh. 

To conclude this analysis, optimization problem \eqref{CapExp2_OF}--\eqref{CapExp2_StrongDuality} is also solved considering a planning horizon of 3 years. The system demand is assumed to increase 5\% every year, and the interest rate is equal to 3\%. Table \ref{table:impactfailure2} provides the optimal investment decisions throughout the planning horizon of two units of 100 MW with a marginal cost of $C^P_{g'}$. The last column of the table provides the relative increase of the producer's profit if equipment failures are accounted for. These results show that considering equipment failures does not only affect the sizing and location of expansion decisions, but also their timing. Observe that as in the target year case discussed in this section, the increase of the GENCO's profit is also significant, reaching 65\% for a marginal cost of 19\$/MWh. 

{\setlength{\tabcolsep}{3mm}\begin{table}[htb!]\begin{center}\renewcommand{\arraystretch}{1}
\caption{Impact of failures on expansion planning (3-year case study)}\label{table:impactfailure2}
\begin{tabular}{c c c c c c c c c c c c}
  \hline 
   $C^P_{g'}$ & \multicolumn{3}{c}{$B^{\text{\tiny NF}}$} & $\overline{\Pi}^{\text{\tiny NF}}$ & \multicolumn{3}{c}{$B^{\text{\tiny F}}$} & $\overline{\Pi}^{\text{\tiny F}}$ & $\Delta\overline{\Pi}$ \\
& $y_1$ & $y_2$ & $y_3$ & & $y_1$ & $y_2$ & $y_3$ & & \\
\hline
15 & $n_2,n_7$ & $-$ & $-$ & 49.97 & $n_1,n_7$ & $-$ & $-$ 	& 53.20 & 6.5\%     \\
16 & $n_2,n_7$ & $-$ & $-$ & 45.64 & $n_1$ & $-$ & $-$ 	& 51.64 & 13.2\%     \\
17 & $n_2,n_7$ & $-$ & $-$ & 41.32 & $n_1$ & $-$ & $-$ 	& 50.79 & 22.9\%     \\
18 & $n_2,n_7$ & $-$ & $-$ & 35.65 & $n_1$ & $-$ & $-$ 	& 49.94 & 40.1\%     \\
19 & $n_7$ & $-$ & $n_1$ & 29.88 & $n_1$ & $-$ & $-$ 	& 49.18 & 64.6\%     \\
20 & $-$ & $n_1$ & $n_7$ & 40.59 & $n_1$ & $-$ & $-$ 	& 48.50 & 19.5\%     \\
21 & $-$ & $n_1$ & $-$ & 46.28 & $n_1$ & $-$ & $-$   	& 47.84 & 3.4\%     \\
22 & $-$ & $-$ & $n_1$ & 33.55 & $n_1$ & $-$ & $-$   	& 47.18 & 40.6\%     \\
23 & $-$ & $-$ & $n_1$ & 33.25 & $n_1$ & $-$ & $-$   	& 46.54 & 40.0\%     \\
24 & $-$ & $-$ & $n_1$ & 33.03 & $n_1$ & $-$ & $-$   	& 46.12 & 39.6\%     \\
  \hline
\end{tabular}\vspace{-3mm}\end{center}
\end{table}
}

%Finally, notice that, for some cases, the construction of the two last units in the table is profitable provided that equipment failures are accounted for.

%It is also worth mentioning that the results provided in Table \ref{table:impactfailure} are significantly influenced by the line capacities of Table \ref{table:linedata}. In the extreme case, if line capacities are such that no network congestion occurs under any scenario, LMPs at all buses are the same and therefore, $\Delta\overline{\Pi}=0$ for all cases.

\vspace{-3mm}
\subsection{Impact of wind speed correlation on expansion planning}\label{sec:ImpWind}

We analyze below the influence of wind speed correlation on the optimal location of new wind farms owned by a GENCO. For the sake of simplicity, equipment failures are not accounted for and just a target year is considered. An existing wind farm of 100MW is located at bus $n_5$. Moreover, new wind farms can only be located at certain nodes of the network, namely, $n_1$, $n_2$, $n_7$, and $n_8$. Wind speed data of 2006 provided by the National Renewable Energy Laboratory (NREL) at five different locations is employed in this analysis. The coordinates of the selected sites are shown in Table \ref{table:coordinateswind}. These data can be freely downloaded from \cite{NREL2010}, and further information on the software employed to simulate these wind power data can be found in \cite{windresorce}. New wind farms will be comprised of 2.5-MW wind generators, model Nordex N80/2500 with a hub height of 105m. The power curve of this turbine model can be found in \cite{danishwind}.

{\setlength{\tabcolsep}{0.5mm}\begin{table}[htb]\begin{center}\renewcommand{\arraystretch}{1}
\caption{Site coordinates for wind data}\label{table:coordinateswind}
\begin{tabular}{c c c c}
  \hline
  $\quad\quad n \quad\quad$    & Coordinates & $\quad\quad n \quad\quad$    & Coordinates \\
  \hline
  $n_1   $    & 42$^{\circ}$21' N, 95$^{\circ}$45' W & $n_7   $    & 42$^{\circ}$26' N, 95$^{\circ}$04' W				\\
  $n_2   $    & 42$^{\circ}$21' N, 95$^{\circ}$20' W & $n_8   $    & 43$^{\circ}$20' N, 95$^{\circ}$18' W					\\
  $n_5   $    & 43$^{\circ}$15' N, 95$^{\circ}$21' W & &					\\
  \hline
\end{tabular}\vspace{-3mm}\end{center}
\end{table}
}

A representative set of 200 wind speed scenarios is used to characterize the correlated wind speed at buses $n_1$, $n_2$, $n_5$, $n_7$ and $n_8$ , being the correlation coefficients those provided in Table \ref{table:corredata}. On the other hand, to isolate the effect of wind correlation on expansion decisions, the marginal distribution of the wind speed at each site should be maintained. For this reason, an uncorrelated scenario set is generated by simply repeating the values of the correlated one as if they were randomly generated, thereby obtaining almost negligible correlation coefficients \cite{feijoo2011simulation}. Similarly to the previous section, the impact of wind correlation on expansion outcomes is determined as follows:

\begin{enumerate}
	\item Model \eqref{CapExp2_OF}--\eqref{CapExp2_StrongDuality} is solved for the uncorrelated scenario set, denoting the wind farm location as $B^{\text{\tiny NC}}$.  
	\item The expected profit if decisions are made disregarding wind speed correlation, $\overline{\Pi}^{\text{\tiny NC}}$, is computed by solving \eqref{CapExp2_OF}--\eqref{CapExp2_StrongDuality} for the correlated scenario set and fixing the wind farm allocation to $B^{\text{\tiny NC}}$. 
	\item The optimal wind farm expansion planning, $B^{\text{\tiny C}}$, and GENCO's expected profit, $\overline{\Pi}^{\text{\tiny C}}$, is calculated by solving model \eqref{CapExp2_OF}--\eqref{CapExp2_StrongDuality} for the correlated scenario set.
	\item $\Delta\overline{\Pi}$, defined as the difference between the profits obtained in 2) and 3), evaluates the impact of wind correlation on expansion decisions. 
\end{enumerate}

{\setlength{\tabcolsep}{0mm}\begin{table}[htb]\begin{center}\renewcommand{\arraystretch}{1}
\caption{Wind speed correlation parameters}\label{table:corredata}
\begin{tabular}{c c c c c c c c c c c }
  \hline
  $n,m$ & $n_1,n_2$ &$n_1,n_5$ &$n_1,n_7$ &$n_1,n_8$ &$n_2,n_5$ &$n_2,n_7$ &$n_2,n_8$ &$n_5,n_7$ &$n_5,n_8$ &$n_7,n_8$ \\
  $\rho_{nm} \quad$ & 0.94 & 0.80 &0.86 &0.81 &0.84 &0.92 &0.84 &0.83 &0.98 &0.83\\
  \hline   
\end{tabular}\vspace{-3mm}\end{center}
\end{table}
}

Table \ref{table:impactcorre} includes the results regarding the optimal location of two wind farms of $N^{T}_{w'}$ wind turbines each. The investment cost of each wind farm is assumed to be proportional to its capacity with a rate of \$1000/kW and independent of its location, being the payback period equal to 40 years. Profits and investment costs are expressed in \$million. The total system load for the target year is approximated by the 20 blocks provided in Table \ref{table:loadblocks}.

{\setlength{\tabcolsep}{3mm}\begin{table}[htb]\begin{center}\renewcommand{\arraystretch}{1}
\caption{Impact of wind correlation on expansion planning (case study)}\label{table:impactcorre}
\begin{tabular}{c c c c c c c}  
  \hline
    $N^{T}_{w'}$ & $C^I_{w'n}$ & $B^{\text{\tiny NC}}$ & $\overline{\Pi}^{\text{\tiny NC}}$ & $B^{\text{\tiny C}}$ & $\overline{\Pi}^{\text{\tiny C}}$ & $\Delta\overline{\Pi}$ \\ 
   \hline     
     100  &12.50 & $n_7,n_8$ & 7.92  & $n_2,n_8$ &  9.21 & 16.3\%  \\ 
     110  &13.75 & $n_7,n_8$ & 5.80  & $n_2,n_8$ &  9.06 & 56.2\%  \\ 
     120  &15.00 & $n_7,n_8$ & 4.07  & $n_2,n_8$ &  8.93 & 119.5\% \\ 
     130  &16.25 & $n_2,n_8$ & 8.87  & $n_8$     &  9.29 & 4.7\%   \\   
  \hline
\end{tabular}\vspace{-3mm}\end{center}
\end{table}
}

Note that if wind speed correlation is neglected, wind farms are located at buses $n_7$ and $n_8$ with comparatively higher LPMs. However, the high wind speed correlation at these two locations ($\rho_{n_7n_8}=0.83$) and the relatively low capacity of the line connecting these two buses ($\overline{P}^F_{n_7n_8}=200$ MW) contributes to high levels of wind spillage for high wind speed scenarios. Consequently, if wind speed correlation is taken into account, this strategy is avoided to reduce wind spillage and increase the profit. Observe that, in some cases, the impact of wind correlation to determine the optimal location of wind farms is of such relevance that GENCO's profit can be increased by 119\%.

\vspace{-3mm}
\subsection{Computational performance}

The simulations presented in this paper are solved using CPLEX 12.3.0 under GAMS on a Linux-based server with eight processors clocking at 1.8 GHz and 20 GB of RAM. The duality gap is set to 0\% in all cases. The approximated computational times required to solve optimization model \eqref{CapExp2_OF}--\eqref{CapExp2_StrongDuality} for the cases presented in sections \ref{sec:ImpFailures} and \ref{sec:ImpWind} are 5-6h and 14-18h, respectively. Note that these computational times can be significantly reduced by applying scenario reduction techniques or  dedicated computational techniques such as decomposition and parallel optimization. However, these issues are out of the scope of this paper and further research is required in this regard.

\section{Conclusions}\label{SectionConclusions}

Driven by the fact that current generation expansion models for GENCOs do not take into account some uncertainties affecting the profitability of an investment, we propose in this paper a bi-level stochastic optimization model to quantify the impact of both equipment failures and wind power production correlation on the expected profit of a GENCO. The proposed model accounts for the variability of the demand throughout the planning horizon, as well as the uncertainty related to both wind speed and unexpected failures of units and transmission lines. Wind speed correlation among different geographical locations is also considered. The use of the primal-dual theory allow us to formulate a single-level mixed-integer equivalent formulation that can be readily solved by off-the-shelf optimization software.

Results provided by the 24-bus case study allows us to identify those cases in which generating expansion decisions can be significantly affected by both equipment failures and wind power correlation, both for a single target year and for the multi-year case. If unit and line failures are modeled, GENCO's expected profit may increase up to 156\%. Likewise, a expected profit increase of 119\% may be incurred if wind speed correlation is taken into account to determine optimal expansion decisions. 

As future research, the proposed model can be solved using a larger test system in order to investigate its scalability features. In this line, the use of dedicated computational methods such as parallel computation or decomposition techniques need to be evaluated. Furthermore, the comparison between the impact of long-term uncertainties such as demand growth or investment costs and the impact of usually disregarded short-term factors such as equipment failures or wind correlation on expansion decisions is also ground for further investigation. Finally, the proposed generation expansion model could be extended so as to incorporate risk measures that allow the power producer to hedge against long-term risks.

\appendix

\section{Linearization}\label{Appendix}

\noindent\emph{- Linearization of $\lambda_{nsby:g\in \Psi_n}P^C_{gsby}$}:

The partial derivative of the Lagrangian function of the lower-level problem with respect to $P^C_{gsby}$ is equal to 0, i.e., 
\begin{align}
& \frac{\partial \Lagr}{\partial P^C_{gsby}} = 0    \Rightarrow \lambda_{nsby:g\in \Psi_n} = C^P_{g} - \phi^{min}_{gsby} - \phi^{max}_{gsby}. \label{Lin1_1} 
\end{align}
The two complementarity conditions corresponding to \eqref{CapExp1_MaxPG} are
\begin{subequations}\begin{align}
&\phi^{min}_{gsby}P^C_{gsby} = 0 \label{Lin1_2a}\\
&(P^C_{gsby}-k_{gs}\overline{P}^C_{g})\phi^{max}_{gsby} = 0. \label{Lin1_2b} 
\end{align}\end{subequations}
Multiplying \eqref{Lin1_1} times $P^C_{gsby}$ and using \eqref{Lin1_2a} and \eqref{Lin1_2b}, we obtain
\begin{equation}
\lambda_{nsby:g\in \Psi_n}P^C_{gsby} = C^P_{g}P^C_{gsby} - k_{gs}\phi^{max}_{gsby}\overline{P}^C_{g}, \label{Lin1_3}
\end{equation}
\noindent which is a linear term.

%\subsection{Linearization of $\lambda_{nsb:w\in\Theta_n}W^S_{wsb}$}
\vspace{2mm}
\noindent\emph{- Linearization of $\lambda_{nsby:w\in\Theta_n}P^W_{wsby}$}:

The partial derivative of the Lagrangian with respect to $P^W_{wsby}$ is equal to 0, i.e.,
\begin{align}
& \frac{\partial \Lagr}{\partial P^W_{wsby}} = 0    \Rightarrow -\lambda_{nsby:w\in\Theta_n} + \gamma^{max}_{wsby} + \gamma^{min}_{wsby} = 0. \label{Lin3_1} 
\end{align}
The complementarity conditions corresponding to \eqref{CapExp1_MaxSpill} are 
\begin{subequations}\begin{align}
&P^W_{wsby}\gamma^{min}_{wsby} = 0 \label{Lin3_2a}\\
&(P^W_{wsby}-\overline{P}^W_{ws})\gamma^{max}_{wsby} = 0. \label{Lin3_2b} 
\end{align}\end{subequations}
Multiplying \eqref{Lin3_1} times $P^W_{wsby}$ and using \eqref{Lin3_2a} and \eqref{Lin3_2b}, we have
\begin{equation}
\lambda_{nsby:w\in\Theta_n}P^W_{wsby} = \gamma^{max}_{wsby}\overline{P}^W_{ws}. \label{Lin3_3}
\end{equation}
%
%
%\subsection{Linearization of $u_{g'n}\lambda_{nsb}P_{g'sb}$}
\vspace{2mm}
\noindent\emph{- Linearization of $\hat{u}_{g'ny}\lambda_{nsby}P^C_{g'sby}$}:

Similarly,
\begin{align}
&\frac{\partial \Lagr}{\partial P^C_{g'sby}} = 0    \Rightarrow - C_{g'} +  \sum_{n}\hat{u}_{g'ny}\lambda_{nsb} = \phi^{min}_{g'sby} - \phi^{max}_{g'sby}. \label{Lin2_1}
\end{align}
The complementarity conditions corresponding to \eqref{CapExp1_MaxPg'} are
\begin{subequations}\begin{align}
&\phi^{min}_{g'sby}P^C_{g'sby} = 0 \label{Lin2_2a}\\
&(P^C_{g'sby}-k_{g's}\sum_n \hat{u}_{g'ny}\overline{P}^{C}_{g'})\phi^{max}_{g'sby} = 0. \label{Lin2_2b} 
\end{align}\end{subequations}
Multiplying \eqref{Lin2_1} times $P^C_{g'sby}$ and using \eqref{Lin2_2a} and \eqref{Lin2_2b}, we arrive at
\begin{equation}
\sum_{n}(u_{g'ny}\lambda_{nsby}-C^P_{g'})P^C_{g'sby} = -\sum_n k_{g's}\hat{u}_{g'ny}\overline{P}^{C}_{g'}\phi^{max}_{g'sby}. \label{Lin2_3}
\end{equation}
%
%\subsection{Linearization of $u_{w'n}\lambda_{nsb}W^{S}_{w'sb}$}
\vspace{2mm}
\noindent \emph{- Linearization of $\hat{u}_{w'ny}\lambda_{nsby}P^W_{w'sby}$}:

The partial derivative of the Lagrangian with respect to $P^W_{w'sby}$ is expressed as
\begin{align}
&\frac{\partial \Lagr}{\partial P^W_{w'sby}} = 0 \Rightarrow   \sum_{n}\hat{u}_{w'ny}\lambda_{nsb} + \gamma^{max}_{w'sby} - \gamma^{min}_{w'sby} = 0. \label{Lin4_1}
\end{align}
The complementarity conditions corresponding to \eqref{CapExp1_MaxSpill2} are
\begin{subequations}\begin{align}
&P^W_{w'sby}\gamma^{min}_{w'sby} = 0 \label{Lin4_2a}\\
&(P^W_{w'sby} - \sum_{n}\hat{u}_{w'ny}\overline{P}^W_{w'ns} )\gamma^{max}_{w'sby} = 0. \label{Lin4_2b} 
\end{align}\end{subequations}
Multiplying \eqref{Lin4_1} times $P^W_{w'sby}$ and using \eqref{Lin4_2a} and \eqref{Lin4_2b}, we obtain
\begin{equation}
\sum_{n}\hat{u}_{w'ny}\lambda_{nsby}P^W_{w'sby} = - \sum_{n}\hat{u}_{w'ny}\gamma^{max}_{w'sby}P^W_{w'sby}. \label{Lin4_3}
\end{equation}
%
%\vspace{2mm}
\noindent \emph{- Linearization of the product of binary and continuous variables:}

Let $\chi$ be a binary variable and $p$ a continuous one bounded by $p^{min}$ and $p^{max}$. Then, the product $z=\chi \cdot p$ is equivalent to the following set of mixed-integer linear expressions:
\begin{subequations}\begin{align}
& z = p - r \label{Lin5_1}\\
& \chi\cdot p^{min} \leq z \leq \chi\cdot p^{max} \label{Lin5_2}\\
& (1-\chi) p^{min} \leq r \leq (1-\chi) p^{max}, \label{Lin5_3}
\end{align}\end{subequations}

\noindent where $r$ is an auxiliary continuous variable.



%% References


\bibliographystyle{model3-num-names}
\bibliography{CapExpansion}






% \end{thebibliography}

\end{document}